\input amstex
\documentstyle{amsppt}
\nologo
\NoBlackBoxes
\NoRunningHeads
\def\p{\hat{p}}
\def\D{\Delta}
\def\dd{\partial\D}
\def\intD{\overset\circ\to D}
\def\Diff{\operatorname{Diff}}
\def\Co{Q}
\pagewidth{28pc}
\pageheight{43pc}
\topskip=\normalbaselineskip

\def\Cl{\mathop{\roman{Cl}}\nolimits}

\def\pt{{\roman{pt}}}
\def\conj{\mathop{\roman{conj}}\nolimits}
\def\oo{\varnothing}
\let\ge\geqslant

\def\C{{\Bbb C}}
\def\R{{\Bbb R}}
\def\Z{{\Bbb Z}}

\def\Cp#1{\C\roman P^{#1}}

\def\Rp#1{\R\roman P^{#1}}
\def\barCP#1{\overline{\C\roman P}^{#1}}
\def\dsum{\bot\!\!\!\bot}
\let\tm\proclaim
\let\endtm\endproclaim
\topmatter
\title
Decomposability of
quotients by complex conjugation for
 rational and Enriques surfaces
\endtitle
\author S. Finashin
\endauthor
\address
\flushpar
Middle East Technical University, Ankara 06531 Turkey
\newline
St.-Petersburg Electrotechnical University,
199376, Russia
\endaddress{}
\email serge\,\@\,rorqual.cc.metu.edu.tr
\endemail
\abstract
The quotients $Y=X/\conj$ by the complex conjugation $\conj\: X\to X$
for complex rational and Enriques surfaces $X$ defined over $\R$
are shown to be diffeomorphic to connected sums of $\barCP2$,
whenever $Y$ are simply connected.
\endabstract
\endtopmatter

\document
\heading
\S1 Introduction
\endheading
\subheading{1.1. The main results}
Given a complex algebraic surface $X$ defined over $\R$, one can consider
{\it the complex conjugation},
$\conj\: X\to X$, and the quotient
$Y=X/\conj$. If $X$ is nonsingular, then $Y$ is a closed 4-manifold,
which inherits an orientation and smooth structure from $X$,
so that the quotient map, $q\: X\to Y$, is 
an orientation preserving and smooth 
double covering branched along the fixed points set, $X_\R$, of
the conjugation.
We call the latter {\it the real part of $X$}.

The main results of this paper are the following theorems.

\tm{1.1.1. Theorem}
If $X$ is a rational nonsingular
complex surface defined over $\R$ with $X_\R\ne\oo$,
then $Y$ is diffeomorphic to $\#_m\barCP2$, 
where $m=\frac12(\chi(X)+\chi(X_\R))-2$.
\endtm

\tm{1.1.2. Theorem}
If $X$ is an Enriques surface defined over $\R$ and $Y$ is simply
connected, then $Y$ is diffeomorphic to $\#_m\barCP2$, 
where $m=\frac12\chi(X_\R)+4$.
\endtm

By making use of the relation between the Euler characteristics
and the signatures of $X$ and $Y$ one gets formulas,
$b_2^+(Y)=p_g$ and $b_2^-(Y)=\frac12(\chi(X)+\chi(X_\R))-2$,
for any algebraic surface $X$ of genus $p_g$ provided $b_1(X)=0$.
Therefore, Freedman's and Donaldson's theorems imply that
$Y$ is {\it homeomorphic} to  $\#_m\barCP2$ if $Y$ is simply connected
and $X$ has genus 0.
Further, $Y$ is simply connected if $X$ is simply connected 
 and $X_\R\ne\oo$.
It can also be simply connected even when $X$ is not,
and the most topological types of real Enriques surfaces
provide such examples.
The goal of the above theorems is, therefore,
 to show that the quotients $Y$ 
for rational and Enriques surfaces $X$ cannot be exotic  $\#_m\barCP2$.

\subheading{1.2. History of the subject}
In the case $X=\Cp2$ the diffeomorphism $Y\cong S^4$ is known as
the Kuiper--Massey theorem \cite{K,M}.
For a real quadric or a cubic surface, $X\subset\Cp3$, the diffeomorphism 
type of $Y$ was determined in \cite{L}.
In \cite{F1,F2} the author found the differential type of $Y$ for several
families of real algebraic surfaces $X$ and, in particular, for certain
rational surfaces.
The Kuiper--Massey theorem implies also that a blow-up
at a point $P\in X_\R$ preserves the quotient, 
$(X\#\barCP2)/\conj\cong Y\,\#\,S^4\cong Y$. 
Blow-up at a pair of conjugated imaginary points of $X$
descend, obviously, to a blow-up of the quotient $Y$. 
It follows that it suffices to prove Theorem 1.1.1 
only for real minimal models of rational surfaces.

\subheading{1.3. Structure of the paper}
In \S2, we recall some basic definitions related to
 real algebraic surfaces and formulate
the classification theorem for real minimal models
of rational surfaces.
In \S3 we analize certain classes of real minimal 
rational surfaces and
complete the proof of Theorem 1.1.1.
In \S4 we deduce Theorem 1.1.2 from Theorem 1.1.1.

\subheading{1.4. Acknowledgements}
I would like to thank A. Degtyarev for the useful discussion concerning
the switch of complex structure arguments \cite{D}
for K3 surfaces. He informed me about \cite{DK2}, where these arguments
are adopted for real Enriques surfaces.

\heading
\S2 Real minimal models for rational surfaces
\endheading
\subheading{2.1. Real surfaces}
By a real variety (a real surface, a real curve)
we mean a pair, $(X,\conj)$, where $X$ is a complex
variety and $\conj\: X\to X$ an anti-holomorphic involution, called
{\it the complex conjugation\/} or {\it the real structure\/}.
A morphism, $X_1\to X_2$,
between real varieties, $(X_i,\conj_i)$, $i=1,2$,
 is called {\it real} 
if it commutes with the complex conjugations in $X_1$ and $X_2$.
It is not difficult to see that if $X$ is an algebraic
complex variety, then there exists a real embedding $X\subset\Cp{N}$,
which makes $X$ an algebraic variety over $\R$
(such an embedding is defined by
$\Cal L\otimes\conj^*(\Cal L)$, where $\Cal L$ 
is a very ample linear bundle on $X$, cf. \cite{CD}).

By a {\it real rational surface} we mean a real surface which is rational
as a complex surface.
A ruled surface $f\: X\to B$ is called {\it real\/} if
$X$ and $B$ are supplied with a real structure so that the ruling $f$ 
is a real morphism.
By a real rational conic bundle we mean
a real morphism $f\: X\to\Cp1$,
which generic fibers are rational curves and 
singular fibers split into wedges of two rational curves.

\subheading{2.2. Real minimal models}
A blow-up $X'\to X$ at a real point of a real surface $(X,\conj)$ 
will be called {\it an elementary real blow-up of type 1.} 
A blow-up at a pair of conjugated imaginary points of $X$
will be called {\it an elementary real blow-up of type 2.}
It is easy to see that in the both cases
$X'$ inherits a real structure making $X'\to X$ a real morphism.
It is well known (cf. \cite{S}) that
any real birational equivalence between 
real surfaces can be decomposed into a sequence of elementary
real blow-ups and downs.
A real surface $(X,\conj)$ is called {\it real minimal} if 
every real birational morphism $X\to X'$, where $X'$ is a smooth
real surface, is an isomorphism. In the other words, $X$ is minimal if
for any exceptional curve, $C\subset X$, 
(i.e. rational curve with $C\circ C=-1$), 
we have $C\circ\conj(C)\ge1$.

As it is mentioned above,
a real blow-up of type 1 does not change the quotient $Y$,
 a real blow-up of type 2 descends to a blow-up on $Y$
and the proof of Theorem 1.1.1 will be completed after we check its
statement for minimal models of real rational surfaces.

According to Comessatti \cite{C}, the minimal models of
real rational surfaces are classified as follows
(see also \cite{CP}, \cite{S}).

\tm{2.2.1. Theorem}
Assume that $X$ is a minimal real rational surface. Then $X$ 
is real isomorphic to one of the following types.
\roster
\item
$\Cp2$ with the usual real structure;
\item
a quadric in $\Cp3$ which does not contain real points
($X_\R=\oo$);
\item
a quadric with $X_\R\cong S^2$;
\item
a real ruled rational surface, $X\cong F_n$, $n\ge 0$,
which has the real part
$X_\R$ homeomorphic to a Klein bottle if $n$ is odd,
and either homeomorphic to a torus or empty if $n$ is even;
\item
a real conic bundle over $\Cp1$ with an even number, $2n$, of 
singular fibers, 
which are all real and consist of pairs of complex conjugated 
exceptional curves;
in this case $X_\R\cong nS^2$;
\item
a real del Pezzo surface of degree 2 with $X_\R\cong 4S^2$,
which can be constructed as a double plane with the
branch locus a real quartic;
\item
a real del Pezzo surface of degree 1 with 
$X_\R\cong \Rp2\dsum 4S^2$, which can be constructed as
a real double quadratic cone, $X\to\Co$,
whose branch locus is the intersections of the cone $\Co\subset\Cp3$,
with a cubic surface,
and a further branch point is the vertex of $\Co$.
\endroster
\endtm

Here $\dsum$ stands for disjoint union and 
$nS^2$ for disjoint union of $n$ spheres.

Note that to prove Theorem 1.1.1 we need to consider only
the cases (4), (5) and (7), since 
in the case (2)  $X_\R=\oo$ and in the cases (1), (3), (6) the statement
of Theorem 1.1.1 is well known (see \cite{K,M} for (1), \cite{L} for (3)
and \cite{F1} for (6)).
The cases (4), (5), (6) are considered in the next section.

\heading
\S3 Proof of Theorem 1.1.1
\endheading

\subheading{3.1. Real ruled surfaces}
Consider a real ruled rational surface 
$p\: X\to \Cp1$. Factorization by $\conj$ gives a smooth map
$\p\: Y\to\D$, where $\D=\Cp1/\conj$, such that the following diagram 
commutes
$$
\CD
X     @>p>>     \Cp1\\
@VqVV              @VVV\\
Y     @>\p>>      \D.\\
\endCD
$$
Split $\D$ into a union of a smaller  disc, $\D_1\subset\D$, and
the annulus, $\Sigma=\Cl(\D-\D_1)$.
It can be easily seen that 
$\p$ is a smooth fibering over the interior of $\D$ with a fiber $S^2$,
therefore $\p^{-1}(\D_1)\cong D^2\times S^2$.
Moreover, $\p^{-1}(\Sigma)\cong S^1\times D^3$.
This is because  the product
of $\p$ and the regular smooth retraction, $\Sigma\to S^1$, 
is a smooth $D^3$-fibering,
 $r\:\p^{-1}(\Sigma)\to S^1$.
Further,
it is well known that a 4-manifold with the boundary $S^1\times S^2$
can be filled up by $S^1\times D^3$ in a unique way,
thus, $Y\cong S^4$.

\subheading{3.2. Conic bundles}
In this section we prove the following

\tm{3.2.1. Proposition}
If $X$ is a real minimal conic bundle with $m$ singular fibers, then
$Y$ is diffeomorphic to $\#_m\barCP2$.
\endtm

\demo{Proof}
Let $p\: X\to\Cp1$ be a real minimal conic bundle and
$\p\:Y\to\D$ obtained by factorization by $\conj$, as in 3.1.
Let $p^{-1}(x_i)=E_i'\cup E_i''$, $i=1,\dots,m$, be the singular fibers
of $X$,  $E_i=p^{-1}(x_i)/\conj=\p^{-1}(x_i)$,
 $E_i'\cap E_i''=\{P_i\}$, $i=1,\dots,m$, 
and $x_1,\dots,x_m$ are ordered
consecutively on the circle $\dd$.
We call $x\in\Rp1\subset\Cp1$ a point of type 1 if
the restriction of $\conj$ to $p^{-1}(x)$ has fixed points 
(hence, $\p^{-1}(x)=p^{-1}(x)/\conj\cong D^2$),
and a point of type 2 if $\conj$ it has no
(then, $\p^{-1}(x)\cong\Rp2$).
It is well known
(and easy to check) that $\p|_{X_\R}$ is a smooth map, which
 has critical points only at $P_i$, $i=1,\dots,m$, and
that the type of a point $x\in\dd$ changes as
we cross $x_i$ (cf. \cite{S}),
which explains why $m$ is even, $m=2n$.

Let $l_i$, $i=1,\dots,m-1$, denote the closed arcs in $\dd$
between $x_i$ and $x_{i+1}$
and  $l_m$ the arc between $x_m$ and $x_1$; all arcs are chosen not to 
contain in their interior the points $x_i$.
Further, let us choose the order of $x_i$ so that $l_i$ consists of
points of type 1 for odd $i$ and type 2 for even $i$.
Consider  small regular closed disjoint neighborhoods
$N_i\subset\D$ of the arcs $l_{2i}$, $i=1,\dots,n$ (see Figure 1),
and put $N_0=\Cl(\D-\cup_{i=1}^nN_i)$,
 $A_i=\p^{-1}(N_i)$, $i=0,\dots,n$.

\tm{3.2.2 Lemma} $A_0$ is diffeomorphic to $S^4-n\intD^4$
(a $4$-sphere with $n$ disjoint open regular $4$-discs removed).
\endtm

\demo{Proof of Lemma $3.2.2$}
Let $T_i$ be  small disjoint closed regular neighborhoods of 
$l_{2i-1}\cap N_0$ in $N_0$, $\D_1=\Cl(N_0-\cup_{i=1}^nT_i)$.
Then  $\p^{-1}(T_i)\cong I\times D^3$
are 3-handles attached to $\p^{-1}(\D_1)\cong D^2\times S^2$,
with the cores ${\pt}\times S^2\subset\d D^2\times S^2$
(see Figure~1).
Such a surgery yields $S^4-n\intD^4$.
\qed
\enddemo

\tm{3.2.3. Lemma} $A_i$, $i=1,\dots,n$, is diffeomorphic to
$\barCP2\#\barCP2-\intD^4$.
\endtm
\demo{Proof of Lemma $3.2.3$}
Let $S_i\subset\D$ be a small closed regular neighborhood of 
$x_i$, $i=1,\dots,2n$ and $B_i=\p^{-1}(S_i)$.
Denote by $H$ the total space of a smooth 2-disc fiber bundle over $S^2$
with the normal number $-2$.
Note, first, that $B_i$ is a regular neighborhood of $E_i$ and 
that $E_i\cong S^2$ has self-intersection $-2$, 
which implies that $B_i\cong H$.
$A_i$ is obtained by gluing together 2 copies of $H$
($B_{2i}$ and $B_{2i+1}$) along the part of their boundary
$\Rp3-\intD^4\subset\Rp3=\d H $.
We can extend the gluing map $\hat f\:\Rp3-\intD^4\to\Rp3-\intD^4$
to $f\:\Rp3\to\Rp3$ and get a closed manifold
$R=H \cup_fH$. Further,  obviously,  $A_i\cong R-\intD^4$ and
since $\Diff_+(\Rp3)$ is known to be connected (cf. \cite{H}),
$Q$ is diffeomorphic to a standard pattern, obtained by
gluing of two copies of $H$ along an orientation revering
diffeomorphism $g\: \Rp3\to\Rp3$
($g$ is orientation reversing since $b_2^+(A_i)=0$).
It is a standard and trivial exercise  in Kirby calculus to check
that if $g$  is a mirror reflection, then
 $H\cup_gH\cong\barCP2\#\barCP2$.
\qed
\enddemo

\midinsert
\topcaption{Figure 1}
Proof of Proposition 3.2.1
\endcaption
\vspace{35mm}
\endinsert

By the above lemmas, $Y$ is obtained from $A_0\cong S^4-n\intD^4$
by filling the ``holes'' with $A_i\cong(\barCP2\#\barCP2)-\intD^4$,
which yields $\#_{m}\barCP2$.
\qed
\enddemo

\subheading{3.3. Del Pezzo surfaces of degree 1}
Figure 2 shows a construction of a curve
$A=\Co\cap C$, the intersection of a quadratic cone $\Co$ with
a cubic surface $C\subset\Cp3$, which
 real part $A_\R$ consists of 4 ovals, i.e. components
contractible in  $\Co_\R-\{P\}$, and a component non-contractible in
in  $\Co_\R-\{P\}$. Consider a double covering $p\: X\to\Co$ branched
along $A$ and at the vertex of $\Co$,
 and choose the one of two possible
liftings to $X$ of the real structures in $\Co$ for which
$p$ maps $X_\R$ into the domain, which consists of
 4 discs bounded by the ovals of $A_\R$ and of
the part of $\Co_\R$ bounded by the non-contractible component and $P$,
as it is shown on Figure~2 (standard details on liftings of real
structures the reader can find, e.g., in \cite{F1}).

\midinsert
\topcaption{Figure 2}
The quadratic cone $\Co$ is the double plane branched along a pair of real 
lines, $L_1$, $L_2$. A cubic surface $C$ is the pull-back of the cubic
curve obtained by a perturbation of 3 lines (left figure).
The real components of $X_\R$ are projected onto the shaded domains
(right figure).
\endcaption
\vspace{50mm}
\endinsert

From the construction of $A$
it can be easily seen that 
a pair of ovals of $A_\R$ can be fused by a deformation of $C$,
after passing a nodal singularity on $A_\R$.
This gives a deformation of double coverings
fusing a pair of real components of $X_\R$.
It is shown in \cite{L} 
that the effect for the quotient $Y=X/\conj$ of such a deformation
is a blow-down (see also \cite{F1}, \cite{W}).
Since Theorem $1.1.1$ is already set up 
for rational surfaces with $4$ real components, 
we have proved it in our case as well.
\qed

\heading 
\S4 Proof of Theorem 1.1.2
\endheading

Let $p\:\widetilde X\to X$ be the double covering of a real Enriques
surface $(X,\conj)$ by K3 surface $\widetilde X$. 
Denote by $c_i\:\widetilde X\to\widetilde X$, $i=1,2$, the anti-holomorphic 
involutions, covering $\conj$;
 $c_1$ and $c_2$
commute and give in product
the covering transformation $t\:\widetilde X\to\widetilde X$ of $p$
 (cf. \cite{DK}).
We assume that one of the involutions, say $c_1$, has a fixed point, since,
otherwise $\Z/2\times\Z/2$ acts freely on $\widetilde X$ and
$Y=X/\conj$ is not simply connected.
Further, $c_2$ (or, equally, $t$) induces an involution
 $c_2'\: X'\to X'$ on the quotient space $X'=X/c_1$.

Now, one can change the complex structure
on $\widetilde X$ so that $c_1$ becomes a holomorphic involution,
whereas $c_2$ and, therefore, $t$ anti-holomorphic.
To get it we pick up a $\conj$-symmetric Kalabi-Yau (Ricci flat)
metric on $\widetilde X$ and following the idea of Donaldson \cite{D}
use that K3 surfaces are hyper K\" ahler to vary the complex
structure on $X$. The only novelty related to Enriques surfaces is that
we need to choose the metric on $X$ to be
symmetric with respect to both $c_2$ and $t$. Such a metric is the pull back
of a conjugation symmetric metric on $X$
(more detailed explanation can be found in \cite{DK2}).

Quotients of K3 surfaces by holomorphic involutions, which have non-empty
fixed point set, are known to be rational (see, e.g., \cite{W}),
therefore, $X'$ gets a structure
of rational surface with the real structure $c_2'$.
Hence, by Theorem 1.1.1, $Y=X'/c_2'$ splits into a connected sum
of $\barCP2$'s provided it is simply connected.
\qed

\enddocument